\documentclass[11pt,a4paper,reqno]{article}
\usepackage{xypic,a4wide,amsmath,amsthm,amssymb}

\newtheorem{thm}{Theorem}[section]

\newtheorem{prop}[thm]{Proposition}
\newtheorem{defin}[thm]{Definition}
\theoremstyle{remark}
\newtheorem{rem}[thm]{Remark}
\newtheorem{ex}[thm]{Example}
\newenvironment{remark}{\begin{rem}\rm}{\qee\end{rem}}
\newenvironment{example}{\begin{ex}\rm}{\qee\end{ex}}
\newcommand{\End}{\mbox{\sl End}}
\newcommand{\diff}{\mbox{\sl Diff}_{\,0}^{\,\,1}}

\newcommand{\rk}{{\mbox{\rm rk}\,}}
\newcommand{\R}{{\mathbb R}}
\newcommand{\C}{{\mathbb C}}

\newcommand{\A}{{\mathfrak A}}
\newcommand{\Q}{{\mathfrak Q}}
\newcommand{\g}{{\mathfrak g}}
\newcommand{\kk}{{\mathfrak k}}
\newcommand{\Lie}{{\mathcal L}}

\newcommand{\sym}{{\mbox{\rm Sym}}}
\newcommand{\vf}{{{\mathfrak X}(M)}}
\newcommand{\qee}{\mbox{\hspace{0.2mm}}\hfill$\triangle$}
\newcommand{\vint}{{\mbox{$\int$}\kern-6pt\mbox{--}}_\epsilon\,}

\begin{document}
\thispagestyle{empty}
\begin{flushright} \small SISSA Preprint 40/2005/fm \\ ITEP-TH-13/05 \\
{\tt arXiv:math.DG/0505392} \end{flushright}
\vfill
\begin{center}
\vfill
{\Large\bf Equivariant cohomology and localization\\[10pt] for Lie algebroids} \\[15pt]
{\sc U. Bruzzo,$^{\ddag\star}$\ L. Cirio,$^{\ddag\star}$\ P. Rossi$^\ddag$} and {\sc V. Rubtsov$^{\P\S}$}\\[5pt]
$^\ddag$ Scuola Internazionale Superiore di Studi Avanzati,\\ Via Beirut 4, 34013 Trieste, Italy\\[3pt] $^\star$
Istituto Nazionale di Fisica Nucleare, Sezione di Trieste\\[5pt]$^\P$ Universit\'e d'Angers,
D\'epartement de Math\'ematiques,\\
UFR Sciences, LAREMA, UMR 6093 du CNRS,\\
2 bd.~Lavoisier, 49045 Angers Cedex 01, France \\[5pt]
$^\S$ ITEP Theoretical Division, \\ 25, Bol.~Tcheremushkinskaya str., 117259, Moscow, Russia
\end{center}
\vfill
\begin{center} \small 16 June  2005, revised  9  October 2006,   4 February 2008 and 28 June 2008 \end{center}
\vfill
\begin{abstract} Let $M$ be a manifold carrying the action of a Lie group $G$, and $A$ a Lie algebroid on
$M$ equipped with a compatible infinitesimal $G$-action. Out of these data we construct an equivariant  Lie
algebroid cohomology and  prove for compact $G$ a related localization formula. As an application we
prove a Bott-type formula.
\end{abstract}
\vfill
\parbox{.75\textwidth}{\hrulefill}\par
{\small \noindent
\par\noindent
U.~B.~and V.~R.~gratefully acknowledge
financial support and hospitality during the respective visits to
Universit\'e d'Angers and SISSA. This work was partially supported by the I.N.F.N.~project PI14.}

 \newpage% \tableofcontents
\section{Introduction}
The notion of Lie algebroid, which may be regarded as a natural generalization of the tangent bundle to
a manifold, allows one to treat in a unified manner
several geometric structures, such as Poisson manifolds, connections on principal bundles, foliations, etc.
Since in many situations one can associate a Lie algebroid to a singular foliation, Lie algebroids seem to provide a
way for generalizing several constructions (e.g., connections) to singular settings. Recently
several field-theoretic models based on Lie algebroids have been proposed, see e.g.~\cite{Str04}.

Every Lie algebroid intrinsically defines a cohomology theory. On
the other hand, the theory of $G$-differential complexes developed
in \cite{Gin99} encompasses the case of the \emph{equivariant Lie
algebroid cohomology} that one can define when a Lie algebroid
carries a (possibly infinitesimal) action of a Lie group $G$ compatible
with an action of $G$ on the base manifold $M$. A natural question
is whether one can generalize to this setting the usual
localization formula for equivariant (de Rham) cohomology. This
seems to have also a certain interest for the applications; for
instance, recently the localization formula has been used to
compute the partition function of $N=2$ super Yang-Mills theory
\cite{Ne03,BFMT03}. In this case the relevant cohomology is the
equivariant de Rham cohomology of the instanton moduli space
--- the moduli space of framed self-dual connections on $\R^4$
(i.e., the relevant Lie algebroid is the tangent bundle to the
instanton moduli space). It seems plausible that super Yang-Mills
theories with a different number of supersymmetry charges may be
treated in a similar way with a different choice of a Lie
algebroid on the instanton moduli space.

In this paper we present a localization formula for equivariant
Lie algebroid cohomology. When a Lie algebroid $A$ on a (compact
oriented) manifold $M$ carries a (infinitesimal) action of a Lie
group, one can introduce an equivariant Lie algebroid cohomology.
If one twists such cohomology by means of the orientation bundle
naturally associated with $A$, equivariant cocycles can be
integrated on $M$. If the group action has only isolated fixed
points, the value of the integral can be calculated as a finite
sum of suitably defined residues at the fixed points.

This localization formula of course reduces to the usual one for
equivariant de Rham cohomology when the Lie algebroid $A$ is the
tangent bundle $TM$. In a similar way, it encompasses a number of
classical localization formulas, providing new proofs for them.
For instance, it implies a generalization to the Lie algebroid
setting of Bott's theorem about zeroes of vector fields (and also  related formulas due to Cenkl and Kubarski \cite{Cenkl73,Kub96}).
This generalizes Bott's formula exactly in the same sense as the notion
of Lie algebroid generalizes that of tangent bundle. Moreover,
when $M$ is complex and $A$ is the Atiyah algebroid of a
holomorphic vector bundle $E$ on $M$, our localization formula can be
specialized so as to produce localization formulas due to Baum-Bott \cite{BB},
Chern \cite{Ch},
Carrell and Lieberman \cite{CL77,Carr78} and K.~Liu \cite{Liu95}.
These aspects will be developed in \cite{BR}. Another localization
formula which is generalized is the one given in \cite{Roub92}.
But many more examples can be given.

The paper is structured as follows. In Section \ref{cohom} we
review the basic definitions and some constructions concerning Lie
algebroid cohomology. Section \ref{equi} introduces the
equivariant Lie algebroid cohomology, basically within the
framework of the theory of \emph{$G$-differential complexes}
developed in \cite{Gin99}. Moreover,  we prove our
localization formula. In Section \ref{connec}, following
\cite{Fern02}, we review one of the several approaches to
connections on Lie algebroids and use this to construct equivariant characteristic classes for Lie algebroids. This will be used in Section
\ref{Bott} to prove a Bott-type formula.

{\bf Acknowledgements.} The authors thank the
referees and the editor A.~A.~Rosly for several useful remarks which helped to substantially improve the presentation and the wording of many basic concepts.
We are grateful to P.~Bressler and Y.~Kosmann-Schwarzbach for their attention to this work.

\par\medskip\section{Lie algebroid cohomology\label{cohom}}
Let $M$ be a smooth manifold. We shall denote
by $\vf$ the space of vector fields on $M$ equipped
with the usual Lie bracket $[\,,]$.
\begin{defin} An algebroid $A$ over $M$ is a vector bundle
on $M$ together with a vector bundle morphism
$a\colon A \to TM$ (called the anchor) and a structure of Lie algebra
on the space of global sections $\Gamma(A)$, such that
\begin{enumerate}
\item $a\colon\Gamma(A)\to \vf $ is a Lie algebra
homomorphism;
\item the following Leibniz rule holds true for every $\alpha$, $\beta\in
\Gamma(A)$ and every function $f$:
$$\{\alpha,f\beta\}=f\{\alpha,\beta\}+a(\alpha)(f)\,\beta$$
(we denote by $\{\,,\}$ the bracket in $\Gamma(A)$).
\end{enumerate}\end{defin}
Morphisms between two Lie algebroids
$(A,a)$ and $(A',a')$ on the same base manifold $M$ are defined in a natural way, i.e.,
they are vector bundle morphisms $\phi\colon A \to A'$ such
that the map   $\phi\colon\Gamma(A)\to\Gamma(A')$ is
a Lie algebra homomorphism, and the obvious diagram
involving the two anchors commutes.

 To any Lie algebroid $A$ one can associate the cohomology complex
$(C^\bullet_A,\delta)$, with $C^\bullet_A=\Gamma(\Lambda^\bullet A^\ast)$
and differential $\delta $ defined by
\begin{multline*}(\delta \xi)(\alpha_1,\dots,\alpha_{p+1}) =
\sum_{i=1}^{p+1}(-1)^{i-1}a(\alpha_i)(\xi(\alpha_1,\dots,\hat\alpha_i,
\dots,\alpha_{p+1})) \\ + \sum_{i<j}(-1)^{i+j}
\xi(\{\alpha_i,\alpha_j\},\dots,\hat\alpha_i,\dots,\hat\alpha_j,\dots,\alpha_{p+1})
\end{multline*}
if $\xi\in C^p_A$ and $\alpha_i \in \Gamma(A), 1\leq i \leq p+1$. The resulting cohomology is denoted by
$H^\bullet(A)$ and is called the cohomology of the Lie algebroid $A$.

\begin{remark} \label{buildalg} One should notice that if $A$ is a vector bundle and $\delta$ is
a derivation of degree $+1$ of the graded algebra  $\Gamma(\Lambda^\bullet A^\ast)$
which satisfies $\delta^2=0$, out of $\delta$
one can construct an anchor $a\colon A \to TM$ and a Lie bracket
on $\Gamma(A)$ making $A$ into a Lie algebroid. One simply defines
$$a(\alpha)(f)=\delta f(\alpha) \qquad\mbox{for}\quad
\alpha\in\Gamma(A),\ f\in C^\infty(M)\,,$$
$$ \xi(\{\alpha,\beta\}) = a(\alpha)(\xi(\beta))-a(\beta)(\xi(\alpha))
-\delta\xi(\alpha,\beta) \qquad\mbox{for}\quad
\alpha,\beta\in\Gamma(A),\ \xi\in \Gamma(A^\ast) \,.$$
\end{remark}

We recall some examples of Lie algebroids.

\begin{example} An involutive distribution inside the tangent bundle
(i.e., a foliation) is a Lie algebroid, whose anchor is injective.
\end{example}
\begin{example} \cite{Va97} Let $\mathfrak M=(M,\mathcal F)$ be a supermanifold; in particular, $\mathcal F$
is a sheaf of $\mathbb Z_2$-graded commutative $\R$-algebras on the differentiable manifold $M$ that can be
realized as the sheaf of sections
of the exterior algebra bundle $\Lambda^\bullet(E)$ for a vector bundle $E$.
Let $D$ be an odd supervector field on $\mathfrak M$ squaring to zero. Then $E^\ast$,
with an anchor and a Lie algebra structure on $\Gamma(\Lambda^\bullet E^\ast)$
given according to Remark \ref{buildalg} by $D$ regarded as a differential for the complex
$\Gamma(\Lambda^\bullet E)$, is a Lie algebroid. Of course, starting from a Lie
algebroid we can construct a supermanifold with  an odd supervector field on  it squaring to zero, so that the
two sets of data are equivalent.
\end{example}
\begin{example} Let $(M,\Pi)$ be a Poisson manifold, where $\Pi$ a Poisson tensor. In this case $A=T^\ast M$ with the Lie bracket of
differential forms
\begin{equation}\label{cdf}
\{\alpha,\beta\}=\Lie_{\Pi(\alpha)}\beta-\Lie_{\Pi(\beta)}(\alpha)
-d\Pi(\alpha,\beta)\end{equation}
(where $\Lie$ is the Lie derivative)
and the anchor is the Poisson tensor.
The cohomology of $A$ is the Lichnerowicz-Poisson cohomology
of $(M,\Pi)$.
\end{example}
\begin{example} Let $P\xrightarrow{p}M$ be a principal bundle with structure
group $G$. One has the Atiyah exact sequence of vector bundles
on $M$
\begin{equation}\label{atialg} 0 \to \mbox{ad}(P) \to TP/G \to TM \to 0\end{equation}
(a connection on $P$ is a splitting of this sequence).
Sections of the vector bundle $TP/G$ are in a one-to-one correspondence with $G$-invariant vector fields on $P$.
On the global sections
of $TP/G$ there is a natural Lie algebra structure, and taking the surjection
$TP/G \to TM $ as anchor map, we obtain a Lie algebroid --- the \emph{Atiyah
algebroid} associated with the principal bundle $P$.

If $E$ is a vector bundle on $M$, we can also associate an Atiyah algebroid
with it: in this case indeed one has an exact sequence
\begin{equation}\label{Atiyah} 0 \to \End(E) \to \diff(E) \to TM \to 0
\end{equation}
where $\diff(E)$ is the bundle of 1-st order differential operators on $E$
having scalar symbol \cite{KS76},
and again $\diff E$, with the natural Lie algebra structure on its global sections and the natural map $\diff(E) \to TM$
as anchor, is a Lie algebroid. The two notions
of Atiyah algebroid coincide when $P$ is the bundle of linear frames of a vector
bundle $E$ (indeed,   an element in $T_uP$ is given by an endomorphism
of  the fibre $E_{p(u)}$ and a vector in $T_{p(u)}M$).

Lie algebroids whose anchor map is surjective, as in the case of the Atiyah
algebroids, are called \emph{transitive}.
\end{example}

Following \cite{ELW99} we now  describe a twisted form of the  Lie algebroid cohomology
together with a natural pairing between the two cohomologies. This will be another
ingredient for the localization formula. Let $Q_A$ be the line bundle $\wedge^r A\otimes\Omega^m_M$, where
$r=\rk A$ and $m=\dim M$, and $\Omega^m_M$ is the bundle of differential $m$-forms on $M$ (later on we shall
denote by $\Omega^m(M)$ the global sections of this bundle). For every $s\in \Gamma(A)$ define a  map
$L_s=\{s,\cdot\}=\colon \Gamma(\wedge^\bullet A)\to
 \Gamma(\wedge^\bullet A)$ by letting
 $$L_s(s_1\wedge \dots\wedge s_k) = \sum_{i=1}^k
 s_1\wedge\dots\wedge\{s,s_i\}\wedge\dots \wedge s_k.$$
Moreover define a map
\begin{eqnarray*} D \colon \Gamma(Q_A) &\to& \Gamma(A^\ast \otimes Q_A)  = \Gamma(A^\ast) \otimes_{C^\infty(M)} \Gamma(Q_A)
\\
D\tau (s) & = & L_s(X) \otimes \mu+ X \otimes \Lie_{a(s)}\mu
\end{eqnarray*}
if $\tau = X\otimes\mu \in \Gamma(Q_A)$ and $s\in\Gamma(A)$.
We consider the twisted complex
$\tilde C^\bullet_A=\Gamma(\wedge^\bullet A^\ast\otimes Q_A) = C^\bullet_A \otimes_{C^\infty(M)}  \Gamma(Q_A)$ with the
differential $\tilde \delta $ defined by
$$\tilde \delta  (\xi\otimes \tau ) = \delta \xi\otimes \tau  + (-1)^{\deg(\xi)}\xi\otimes D\tau$$ and
$\xi \in C^\bullet_A = \oplus_{k=0}^{r}\Gamma(\wedge^{k}(A^\ast)\,.$
We denote the resulting cohomology $H^\bullet(A,Q_A)$.

There is a naturally defined map\footnote{We use interchangeably the notations $\psi\lrcorner$ and $i_\psi$ to denote  the inner product with an element $\psi$, according to typographical convenience.}
  $p\colon  \tilde C^\bullet_A \to \Omega^{\bullet-r+m}(M)$
$$p(\psi\otimes X \otimes \mu) =
(a(\psi\lrcorner X))\lrcorner \mu\,.$$
\begin{prop} \label{pischain} The morphism $p$ is a chain map,\footnote{We are thankful to A. Rosly for pointing out this fact, and for calling to our attention the following proof.}
 in the sense that the diagram
\begin{equation}\label{pseudochain}  \xymatrix{
 \tilde C^{k}_A  \ar[r]^{p\ \ \ \ }\ar[d]^{\tilde\delta} & \Omega^{k-r+m}(M) \ar@<-7mm>[d]_{d} \\
  \tilde C^{k+1}_A  \ar[r]^{p\ \ \ \  } & \Omega^{k-r+m+1}(M)
}\end{equation}
commutes up to a sign, i.e., on $C^{k}_A$ one has
\begin{equation}\label{chainmap} p\circ\tilde\delta=(-1)^k\  d\circ p\,.\end{equation}
\end{prop}
\begin{proof}
If we define a Lie derivative on $\tilde C^\bullet_A$ by letting
$L_s=i_s\circ\tilde\delta+\tilde\delta\circ i_s$ for elements $s\in\Gamma (A)$, this satisfies
the commutation relation
 \begin{equation}\label{commLie} p\circ L_s=\Lie_{a(s)}\circ p \end{equation}
 on $\tilde C^k_A$. By using this identity we may prove \eqref{chainmap} by descending induction
 on $k$. For $k=r$ the identity reduces to $0=0$. For smaller values of $k$ it is enough to prove
 the identity for those $c'\in  \tilde C^k_A$ that can be represented as $c'=i_s\,c$ for
 some $s\in\Gamma(A)$ and $c\in  \tilde C^{k+1}_A$. In this case the result follows from
 a simple computation which uses \eqref{commLie}.
\end{proof}

Let $M$ be compact and oriented, and note that
 since   $\tilde C^r_A\simeq\Omega^m(M)$ (canonically),
one can integrate elements of $\tilde C^r_A$ over $M$.
There is a nondegenerate pairing
\begin{eqnarray*} C^k_A\otimes_{C^\infty(M)}  {\tilde C}^{r-k}_A &\to& \R \\
\xi \otimes (\psi\otimes X \otimes \mu) & \mapsto &   \int_M(\xi\wedge\psi,X)\,\mu\,.
\end{eqnarray*}
A version of Stokes' theorem holds  for the
complex $\tilde C_A^\bullet$  \cite{ELW99}:
if  $c\in\tilde C^{r-1}_A$, then
$$\int_M\tilde \delta c=0\,.$$
This formula follows from the identity \eqref{chainmap} for $k=r-1$.
It implies that the pairing descends to cohomology, yielding a bilinear map
\begin{equation}\label{pairing} H^\bullet(A)\otimes H^{r-\bullet}(A,Q_A)\to\R\,.
\end{equation}
This pairing in general may be degenerate.

One also has a natural morphism $C^\bullet_A\otimes_{C^\infty(M)} \tilde C^\bullet _A \to
\tilde C^\bullet _A $ which is compatible with the degrees. Again this descends
to cohomology and provides a cup product
\begin{equation}\label{cup} H^i(A)\otimes H^j(A,Q_A)\to H^{i+j}(A,Q_A)\,.
\end{equation}

\par\medskip\section{Equivariant  cohomology and localization\label{equi}\label{loca}}
In this section we introduce an equivariant cohomology for Lie algebroids,
basically following the pattern exploited in \cite{Gin99} to define
equivariant cohomology for Poisson manifolds (and falling within the
general theory of equivariant cohomology for $G$-differential complexes
developed there). Moreover, we write a localization formula for the equivariant
Lie algebroid cohomology (Theorem \ref{loc}).

We assume that there is an action of a Lie algebra $\g$ on the Lie algebroid $A$, i.e., that there is a Lie algebra map
\begin{equation}  b\colon\g\to\Gamma(A)\,. \end{equation}
By composing this with the anchor map
we obtain a Lie algebra homorphism
 $\tilde\rho=a\circ b \colon\g\to \vf $, i.e.,
 an action of $\g$ on $M$. Lie algebra maps like our $b$ have been introduced in \cite{Mac05} in the case
of the Atiyah algebroids under the name of ``derivation representations.''

\begin{example}\label{exval} (Cf.~ \cite{Vai90})  Let $\Pi$ be a regular Poisson tensor on $M$ (i.e., $\Pi$ has constant rank), and let $\mathcal S=\mbox{Im}(\Pi)$
be the associated symplectic distribution. The family of
symplectic forms defined on the leaves of the distribution $\mathcal
S$ yields an isomorphism $\mathcal S\simeq \mathcal S^\ast$. Then $\mathcal S^\ast$ is
a subbundle of $TM$ and   $\Gamma(\ker\Pi)$ is
an ideal in $\Omega^1(M)$ with respect to the Lie algebroid structure in $T^\ast M$ given by the bracket in Eq.~\eqref{cdf}. Moreover,  $\mathcal S^\ast$ is a Lie subalgebroid of $TM$;
its cohomology is called
the   \emph{tangential Lichnerowicz-Poisson
cohomology}.  Assume
now that $M$ carries the action $\rho$ of a Lie group $G$;
if $\g$ is the Lie algebra of $G$, for every $\xi\in \g$  let us denote by
$$\xi^\ast = \frac{d}{dt} \rho_{\exp(-t\xi)\vert t=0}$$
the corresponding fundamental vector field (thus, we have a Lie
algebra homomorphism $\tilde \rho\colon \g\to\vf$, $\tilde\rho(\xi)=\xi^\ast$).
If the  $G$-action   is tangent to $\mathcal
S^\ast$,   we obtain an infinitesimal action $b\colon\g\to\Gamma(\mathcal S^\ast)$
(this is what   in \cite{Gin99} is called a \emph{cotangential
lift}).
\end{example}

If for a $\xi\in\g$ the point $x\in M$ is a zero of $\xi^\ast$ then we have the usual endomorphism
$$ L_{\xi,x}  \colon T_pM \to T_pM,\qquad L_{\xi,x}  (v)=[\xi^\ast,v] \,.$$

We consider the graded vector space
$$\A^\bullet=\mbox{Sym}^\bullet(\g^\ast)\otimes\Gamma(\wedge^\bullet A^\ast)$$
with the grading
$$\deg({\mathcal P}\otimes\beta) = 2\deg({\mathcal P})+\deg\beta$$
where ${\mathcal P}\in \mbox{Sym}^\bullet(\g^\ast)$ and $\beta\in \Gamma(\wedge^\bullet A^\ast) \,.$

We will consider $\mathcal P$ as a polynomial function on $\g$ and define the equivariant differential $\delta_\g\colon \A^\bullet\to \A^{\bullet+1}$
\begin{equation}\label{eqdiff} (\delta_\g({\mathcal P}\otimes\beta))(\xi) = {\mathcal P}(\xi)\left( \delta (\beta) - i_{b(\xi)} \beta\right),
\end{equation}
where both sides have been evaluated on an element $\xi\in\g$.
If we denote   $\A_G^\bullet = \ker \delta_\g^2$, then
$(\A_G^\bullet,\delta_\g)$ is a   complex, whose cohomology
we denote $H^\bullet_G(A)$ and call the \emph{equivariant cohomology}
of the Lie algebroid $A$ (to be more precise, of the pair $(A,b)$).

By considering the graded vector space
$$\Q^\bullet = \A^\bullet\otimes \Gamma(Q_A) =  \mbox{Sym}^\bullet(\g^\ast)\otimes\Gamma(\wedge^\bullet A^\ast\otimes Q_A)$$
with a differential $\tilde \delta_\g$ obtained by coupling $\delta_\g$ with the
differential $D$, and letting $\Q^\bullet_G=\ker \tilde \delta_\g^2$,
one also has a twisted equivariant cohomology
$H^\bullet_G(Q_A)$, and there is a cup product
$$H^i_G(A)\otimes H^k_G(Q_A) \to H^{i+k}_G(Q_A).$$

We write now a localization formula. In view of Proposition \ref{pischain}, its right-hand
side   can be computed by means of the usual localization formula in equivariant de Rham cohomology;
the integral of an equivariantly closed  $\Q^\bullet_G$-cocycle $\gamma(\xi)$ is actually the integral of the differential
form $p(\gamma(\xi))$, and $p(\gamma)$ is a cocycle
in the equivariant de Rham complex due to Proposition \ref{pischain}.  It is indeed
quite easy to prove the identity
\begin{equation}\label{identity}
p(\tilde\delta_\g(\gamma)) = (-1)^k d_\g (p(\gamma))
\end{equation}
when $\gamma\in\Q^k_G$.  Here
 $d_\g $ is the differential in the usual equivariant de Rham cohomology. This follows
from Proposition \ref{pischain} and the equalities
\begin{eqnarray*} i_{\xi^\ast}p(\gamma) &=& i_{\xi^\ast}\left[a(\psi\lrcorner X)\lrcorner\mu\right]
= (\xi^\ast\wedge a(\psi\lrcorner X))\lrcorner\mu \\
&=& a(b(\xi)\wedge(\psi\lrcorner X))\lrcorner\mu = (-1)^{k-1}a\left((i_{b(\xi)}\psi)\lrcorner X\right)\lrcorner\mu \\
&=& (-1)^{k-1}p(i_{b(\xi)}\gamma)
\end{eqnarray*}
having set $\gamma=\psi\otimes X\otimes\mu$.

 Let $M$ be a closed manifold which carries
the action $\rho$ of a compact Lie group $G$.
We also assume that $M$ is oriented, and that a $\xi\in\g$ has been
chosen such that $\xi^\ast=\tilde\rho(\xi)$ (this has been defined in Example \ref{exval}) only has isolated zeroes.  We denote by $M_\xi$ the set of
such zeroes. Note that due to the compactness of $G$, we have $\det(L_{\xi,x})\ne 0$ at every
isolated zero $x$, and the dimension $m$ of $M$ is necessarily even (as we shall assume henceforth).

If the rank $r$ of $A$ is smaller than the dimension $m$ of $M$, then for every equivariantly closed
 $\gamma\in\mathfrak Q^\bullet$ we have $\int_M\gamma(\xi)=0$ for dimensional reasons, since
 $p(\gamma)_0=0$ in that case.  (Here the subscript $0$ denotes the piece of degree 0 in the usual de Rham grading).
 We may therefore henceforth assume that $r\ge m$.

\begin{thm}\label{loc} Let $M$ be a closed oriented $m$-dimensional manifold on which  a compact Lie group $G$ acts.  Let $A$ be a rank $r$ Lie algebroid on $M$, with $r\ge m$, and assume that a Lie algebra homomorphism
$b\colon \g \to\Gamma(A)$ exists making the diagram
\begin{equation}\label{infeq} \xymatrix{\g\ar[r]^b \ar[dr]_{\tilde\rho}& \Gamma(A) \ar[d]^a \\ & \vf }\end{equation}
 commutative; in other terms, $\tilde\rho$ is the Lie algebra homomorphism
 $\tilde\rho(\xi)=\xi^\ast$ (here $\g$ is the Lie algebra of $G$). Moreover, assume that $\xi\in \g$
is such that the associated fundamental vector field $\xi^\ast$ has only isolated zeroes. Finally, let $\gamma\in\mathfrak Q^\bullet$
be equivariantly closed, $\tilde\delta_\g\gamma=0$.

Then  the following localization formula holds:
\begin{equation}\label{locf}\int_M\gamma(\xi) = (-2\pi)^{m/2} \sum_{x\in M_\xi}
\frac{p(\gamma(\xi))_0(x)}{\operatorname{det}^{1/2} L_{\xi,x}}
.\end{equation}
\end{thm}
\begin{proof}
Since in the left-hand side we are actually integrating the differential form $p(\gamma(\xi))$,
the formula follows from the identity \eqref{identity} and the usual localization formula.
\end{proof}

One easily checks that in the case of the ``trivial'' algebroid
given by the tangent bundle with the identity map as anchor,
 this reduces to the ordinary localization formula for the equivariant de Rham
cohomology (see e.g.~\cite{BGV92}).

\begin{remark} \label{rank} If $r\ge m$ and the rank of the linear morphism $a$ at the point $p$
is not maximal (i.e., less than $m$), then $p(\gamma(\xi))_0(x)=0$.
\end{remark}

\begin{remark} \label{vectfield} As a particular case of Theorem \ref{loc}, one can  state a localization formula related to the action
of a vector field on $M$.
Let $M$ be a compact oriented $m$-dimensional manifold,  and $X \in \Gamma(\mathrm{T}M)$  a vector field on $M$
with isolated zeroes which generates a circle action.
Let $A$ be a rank $r$ Lie algebroid on $M$ such that there exists $\tilde{X}\in \Gamma(A)$ with $a(\tilde{X})=X$
(where $a$ is the anchor map).
Then, for the integration of a form
$\gamma \in \Gamma(\wedge^\bullet A^*\otimes \wedge^r A\otimes \wedge^m \mathrm{T}^*M)$
such that $\tilde{\delta}_X \gamma:=(\tilde{\delta}-i_{\tilde{X}})\gamma=0$,
the localization formula (\ref{loc}) holds.

One can replace the assumption that $X$ generates a circle action
by assuming that $X$ is an isometry of a Riemannian metric on $M$.
\end{remark}

\par\medskip\section{Connections and characteristic classes for Lie algebroids\label{connec}}
Several applications of the localization formula may be given by using
a notion of a characteristic class for a Lie algebroid. We start by introducing the concept of \emph{$A$-connection,}
cf. \cite{Fern02}.

Let $A$ be a Lie algebroid with anchor $a$, and let
$P\xrightarrow{p}M$ be a principal bundle with structure group
$K$. Note that the pullback $p^\ast A = A \times_M P$ carries a
natural $K$-action, and $A \simeq p^\ast A/K$. Also the tangent
bundle $TP$ carries a natural $K$-action, and one has
$p_\ast(vk)=p_\ast(v)$ for $v\in TP$ and $k\in K$, so that one
has an
induced map $p_\ast\colon TP/K \to TM$  which is the anchor of the Atiyah
algebroid associated with $P$, see Eq.~(\ref{atialg}).

\begin{defin} An $A$-connection on $P$ is a bundle map $\eta\colon p^\ast A\to TP$ such that:
\begin{enumerate} \item  the following diagram commutes:
$$\xymatrix{ p^\ast A \ar[d] \ar[r]^\eta & TP \ar[d]^{p_\ast} \\ A \ar[r]^{a} & TM}$$
\item $\eta$ is $K$-equivariant, i.e., $\eta(uk,\alpha) =
R_{k\ast}\eta(u,\alpha)$ for all $k\in K$, $u\in P$, $\alpha\in
A$. (Here $R_k$ denotes the structural right action of an element
$k\in K$ on $P$.) \end{enumerate} If $P$ is the bundle of linear
frames of $A $, $\eta$ is called an $A$-linear
connection.\end{defin}

By its equivariance, an $A$-connection $\eta$ for a principal $K$-bundle $P$  defines a bundle map
$\omega_\eta\colon A \to TP/K$, called the \emph{connection 1-section} of
$\eta$.  One has $p_\ast\circ \omega_\eta = a$.

\begin{remark} \label{remcon} 1. The usual notion of connection is recovered by
taking $A=TM$, and $\eta$ is then the corresponding horizontal lift $\eta\colon
p^\ast TM\to TP$.

2. An ordinary connection   on  $P$ (regarded as the associated
horizontal lift $\zeta\colon p^\ast TM\to TP$) defines an
$A$-connection $\eta$ for $P$  by letting $\eta = \zeta\circ
p^\ast a$. \end{remark}

If $E\xrightarrow{p_E} M$ is a vector bundle associated with $P$
via a representation of $K$ on a linear space, an $A$-connection
on $P$ defines a similar structure on $E$, that is, a bundle map
$\eta_E\colon p_E^\ast A\to TE$ which makes the diagram
$$\xymatrix{ p^\ast_E A \ar[d]_A \ar[r]^{\eta_E} & TE \ar[d]^{p_{E\ast}} \\ A \ar[r]^{a} & TM}$$
commutative. The $A$-connection $\eta_E$ defines in the usual way
a covariant derivative $\nabla\colon\Gamma(A
)\otimes_\R\Gamma(E)\to\Gamma(E)$: if $\phi\colon TP/K\to\diff(E)$
is the natural map,\footnote{A section $X$ of $TP/K$ is a
$K$-invariant vector field on $P$. Since there is an obvious  map
$TP\to\diff(P\times V)$, where $V$ is the standard fibre of $E$,
by equivariance $X$ yields a differential operator on $E$ whose
symbol is scalar.} one sets
$\nabla_\alpha=(\phi\circ\omega_\eta)(\alpha)$. This covariant
derivative satisfies the Leibniz rule
\begin{equation}\label{leibniz}
\nabla_\alpha(fs) = f\nabla_\alpha(s) + a(\alpha)(f)s
\end{equation}
for all functions $f$ on $M$.

Let us introduce the notion of a \emph{$G$-equivariant} $A$-connection.
Assume that a Lie group $G$ acts on $M$, and that this action $\rho$ lifts
to an action $\hat \rho$ on $A$; this means that for every $g\in G$ we
have a commutative diagram
$$\xymatrix{ A \ar[d]_a \ar[r]^{\hat\rho_g} & A  \ar[d]^a \\
TM \ar[r]^{\rho_{g\ast}} & TM}$$
Moreover, we assume that $\rho$ also lifts to an action $\tilde\rho$ on the principal
$K$-bundle $P$ which commutes with the structural $K$-action.

\begin{defin}\label{invconn}
A $G$-equivariant $A$-connection $\eta$ on $P$ is an $A$-connection $\eta$
such that the diagram
$$\xymatrix
{ p^\ast A \ar[r]^{\hat\rho_g} \ar[d]_\eta & p^\ast A \ar[d]^\eta \\
TP \ar[r]^{\tilde\rho_{g\ast}} & TP }$$
commutes for every $g\in G$.
\end{defin}
Since the action of $G$ on $P$ commutes with the action of $K$,
we have an induced action $\tilde\rho_\ast$ on $TP/K$, and the condition
for $\eta$ to be $G$-equivariant may be stated in terms of the connection 1-section
$\omega_\eta$ as the commutativity of the diagram
$$ \xymatrix { A \ar[r]^{\hat\rho_g} \ar[d]_{\omega_\eta} & A \ar[d]^{\omega_\eta} \\
TP/K \ar[r]^{\tilde\rho_{g\ast}} & TP/K}$$
%If one only has an infinitesimal action of $G$ on $A$ (in the sense of
%Theorem \ref{loc}) then the $G$-equivariance of an $A$-connection $\eta$
%amounts to requiring that its connection 1-section $\omega_\eta$ satisfies
%the condition
%$$ \omega_\eta(\{b(\xi),\alpha\}) = [\tilde\xi^\ast,\omega_\eta(\alpha)]$$
%for all $\xi\in\g$ and $\alpha\in\Gamma(A)$, where $\tilde\xi^\ast$
%is the section of $TP/K$ induced by the vector field on $P$ which
%generates the action of $G$.

The \emph{curvature} ${\mathcal R}_\eta$ of an $A$-connection $\eta$
on a principal $K$-bundle $P$ may be defined in terms of the map
$\omega_\eta$ as the element in $\Gamma(\Lambda^2A^\ast\otimes
TP/K)$ given by $${\mathcal R}_\eta(\alpha,\beta) =
[\omega_\eta(\alpha),\omega_\eta(\beta)] -
\omega_\eta(\{\alpha,\beta\}).$$ As a matter of fact $p_\ast\circ
{\mathcal R}_\eta=0$, so that ${\mathcal R}_\eta$ is an element in
$\Gamma(\Lambda^2A^\ast\otimes \mbox{ad}(P))$. The curvature $\mathcal R_\eta$
satisfies analogues of the \emph{structure equations} and \emph{Bianchi} identities.
These identities are conveniently stated in terms of the
so-called \emph{exterior $A$-derivative}
\begin{eqnarray*}
D_A & \colon & \Gamma(\Lambda^\bullet A^\ast\otimes TP/K) \to
\Gamma(\Lambda^{\bullet +1} A^\ast\otimes TP/K) \\
(D_A\chi )(\alpha_1,\dots,\alpha_{p+1})& =&
\sum_{i=1}^{p+1}(-1)^{i-1} [\omega_\eta(\alpha_i),\chi(\alpha_1\dots,\hat\alpha_i,
\dots,\alpha_{p+1})] \\ &+&  \sum_{i<j}(-1)^{i+j}
\chi(\{\alpha_i,\alpha_j\},\dots,\hat\alpha_i,\dots,\hat\alpha_j,\dots,\alpha_{p+1})
\end{eqnarray*}
in the form
\begin{equation}\label{ident} \mathcal R_\eta= D_A\omega_\eta - \tfrac12[\omega_\eta,\omega_\eta]\,,
\qquad D_A\mathcal R_\eta=0\,.\end{equation}
Note that the equivariance of the connection may be expressed by the condition
$$[D_A,\Lie_{\tilde\xi^\ast}]=0$$
for all $\xi\in\g$. Here $\Lie_{\tilde\xi^\ast}$ is the Lie derivative of sections of $TP/K$
with respect to the section $\tilde\xi^\ast$
of $TP/K$ induced by the vector  field on $P$ which
generates the action of $G$ (i.e., $ \Lie_{\tilde\xi^\ast}(v)=[\tilde\xi^\ast,v]$).

If the connection $\eta$ is equivariant, we may equivariantize these relations, in particular
by defining the equivariant exterior $A$-derivative $D_A^{\g}$, which acts on $\operatorname{Sym}^\bullet(\g^\ast)\otimes
 \Gamma(\Lambda^\bullet A^\ast\otimes TP/K)$ as
 $$(D_A^{\g}\chi)(\xi) = D_A(\chi(\xi))-(i_{\tilde\xi^\ast}\otimes\operatorname{id})(\chi(\xi))\,.$$
 Moreover we define the equivariant curvature $R_\eta^{\g}$ of $\eta$ as
 $$(R_\eta^{\g} \chi)(\xi)= R_\eta (\chi(\xi))+ \Lie_{\tilde\xi^\ast}(\chi(\xi)) - [D_A,i_{\tilde\xi^\ast}\otimes\operatorname{id}]
  (\chi(\xi)) = R_\eta (\chi(\xi)) + \mu(\chi(\xi))$$
  where the last equality defines the ``moment map'' $\mu$. Moreover, the square brackets in this equation denote an anticommutator. An easy calculation shows that the equivariant curvature satisfies
  the equivariant Bianchi identity
\begin{equation}\label{equivBianchi}  D_A^{\g} R_\eta^{\g} = 0\,.\end{equation}

We can also write  the identities \eqref{ident} and \eqref{equivBianchi} in a local form, involving
the \emph{local connection 1-sections} defined in the following way.
Let $\{U_i\}$ be an open cover of $M$ over which the bundle $P$
trivializes. Then one has local isomorphisms
$$\psi_j\colon (TP/K)_{\vert U_j} \to TU_j\times\mathfrak k\,,$$
where $\mathfrak k$ is the Lie algebra of $K$.
The local connection 1-sections are defined by the condition
$$\omega_j (\alpha ) = \mbox{pr}_2\circ \psi_j \circ \omega_\eta(\alpha)$$
(where $ \mbox{pr}_2$ is the projection onto the second factor
of $ TU_j\times\mathfrak k$),
and one analogously defines the local curvature 2-sections $\mathcal R_j$.
The identities \eqref{ident} are now written in the form
\begin{equation}\label{identloc} \mathcal R_j=\delta\omega_j +\tfrac12[\omega_j,\omega_j]\,,
\qquad \delta \mathcal R_j + [\omega_j,\mathcal R_j]=0\,.
\end{equation}

In the same way, the equivariant curvature may be represented
 by local 2-sections ${\mathcal R}_j^{\g}$ which,
in view of Eq.~\eqref{identloc}, satisfy the    identities
\begin{equation} \label{equivcurv}
\qquad \delta_{\g} {\mathcal R}_j^{\g}  + [\omega_j,{\mathcal R}_j^{\g} ] =0\,.
\end{equation}

The Chern-Weil homomorphism is defined as follows. Let
$I^\bullet(\kk)=(\sym^\bullet \kk^\ast)^K$ ,
choose an $A$-connection $\eta$ for $P$,  and for any polynomial $Q\in I^\ell(\kk)$  of degree $\ell$
 define the
element $\lambda_Q\in C_A^{2\ell}$
$$\lambda_Q(\alpha_1,\dots,\alpha_{2\ell}) =
\sum_{\sigma}(-1)^\sigma \tilde Q({\mathcal
R}_\eta(\alpha_{\sigma_1},\alpha_{\sigma_2}), \dots, {\mathcal
R}_\eta(\alpha_{\sigma_{2\ell-1}},\alpha_{\sigma_{2\ell}}))
$$
where the summation runs over the permutations of $2\ell$ objects and $\tilde Q$
is the polarization of $Q$, i.e., the unique Ad-invariant symmetric function of $\ell$
variables in $\kk$ such that
$ \tilde Q(\chi ,\dots, \chi) = Q(\chi )$ for all $\chi\in\kk$.
 One
proves  that this cochain is $\delta$-closed, and that the resulting
cohomology class $[\lambda_Q]$ does not depend on the connection,
thus defining a graded ring homomorphism $\lambda\colon I^\bullet(\kk) \to H^{2\bullet} (A)$.
If $\tilde\lambda\colon  I^\bullet(\kk) \to H^{2\bullet}_{dR}(M)$
is the usual Chern-Weil homomorphism to the de Rham cohomology of $M$,
then there is a commutative diagram
\begin{equation}\label{cwhs}\xymatrix{  I^\bullet(\kk) \ar[r]^{\tilde\lambda}\ar[rd]_{\lambda} &
H^{2\bullet}_{dR}(M) \ar[d]^{a^\ast} \\ & H^{2\bullet }(A) }\end{equation}
Using the Chern-Weil homomorphism one can introduce various
sorts of characteristic classes for the Lie algebroid $A$. However
due to diagram (\ref{cwhs}), and somehow unpleasantly, these
characteristic classes are nothing more that the image under $a^\ast$
of the corresponding characteristic classes for the bundle $E$. To show this,
choose any (ordinary)  connection on $P$ to compute a characteristic
classes in $H^\bullet_{dR}(M)$, and the $A$-connection associated with
it (see Remark \ref{remcon} (2)) to compute a characteristic class
in $H^\bullet (A) $. The two characteristic classes are then related
by the morphism $a^\ast$ as in diagram \eqref{cwhs}.

In the following we shall use ``Pontryagin-like'' characteristic classes:
namely, we take   $\kk=\mbox{gl}(r, \C)$,
so that $P$ is the bundle of linear frames of a complex vector bundle
$E$. We assume that $E$ is the complexification of a real vector bundle.
Let   $Q_i$ be the $i$-th elementary
Ad-invariant polynomial, and denote by $\lambda_i$
the corresponding characteristic class. These characteristic classes vanish when $i$ is odd (to check this, take for instance a connection on $P$ which is compatible with a fibre
metric on $E$, and compute the characteristic classes by means of the
induced $A$-connection).

If $Q\in I^\bullet (\mbox{gl}(r, \C))$ is  an  Ad-invariant homogeneous symmetric
polynomial of degree $2i$ on the Lie algebra $\mbox{gl}(r, \C)$,  one has
$Q({\mathcal R}^{\g}_\eta)\in \A^{4i}_G$. The following statement is easily
proved.
\begin{prop} The element $Q({\mathcal R}^{\g}_\eta )$ is $\delta_\g$-closed.
The corresponding cohomology class $\lambda_Q^{\g}(A) \in H^{4i}_G(A)$
does not depend on the equivariant connection $\eta$.
\end{prop}
\begin{proof}  First one shows that
the quantity $Q({\mathcal R}^{\g}_\eta )$ is $\delta_{\g}$-closed by using the equations \eqref{equivcurv}. To prove the second claim,
if  $\eta$, $\eta'$ are two equivariant $A$-connections for the principal bundle $P$, define the 1-parameter family of connections
$$\eta_t = t\,\eta' + (1-t)\,\eta$$
with $0 \le t \le 1$, and let
\begin{multline*} q(\eta,\eta') (\alpha_1,\dots,\alpha_{2i})=  \\ i  \sum_{\sigma}(-1)^\sigma  \int_0^1 \left[\tilde Q(\frac{d}{dt}\omega_{\eta_t}\left(\alpha_{\sigma_1},\alpha_{\sigma_2}),
{\mathcal R}^{\g}_{\eta_t}(\alpha_{\sigma_{3}},\alpha_{\sigma_{4}}),\dots
{\mathcal R}^{\g}_{\eta_t}(\alpha_{\sigma_{2i -1}},\alpha_{\sigma_{2i }})\right)\right]\,dt
\end{multline*}
where $\tilde Q$ is the polarization of the polynomial $Q$.
A straightforward computation which uses    the identity \eqref{equivcurv}
and the identity
$$ \frac{d}{dt}{ \mathcal R}^{\g}_{\eta_t} = D_A^{\g}  \frac{d}{dt} \omega_{\eta_t}$$
now shows that
$$ Q({\mathcal R}^{\g}_{\eta'} ) - Q({\mathcal R}^{\g}_\eta ) = \delta_{\g} q(\eta,\eta')$$
whence the claim follows.
  \end{proof}

When $Q$ is the standard $i$-th elementary Ad-invariant polynomial
$\varsigma_i$ on  $\mbox{gl}(r, \C)$ we shall denote by
$\lambda_i^{\g}(A)$ $(i=1,\dots,r)$ the corresponding equivariant
characteristic class (vanishing for odd $i$). As we discussed previously, these
are just the images under the morphism $a^\ast$
(the adjoint of the anchor map) of the equivariant Chern classes of the (complexified) vector
bundle $E$.

\par\medskip\section{A Bott-type formula\label{Bott}}
As an application of our localization theorem we prove a result which
generalizes the classical Bott formula \cite{Bott67} as well as similar
results by Cenkl and Kubarski \cite{Cenkl73,Kub96}. Bott's formula comes
in different flavours, according to the assumptions that one makes on the
vector field entering the formula. The case we consider here extends
the usual Bott formula for a vector field which generates a circle action and has isolated critical points.

Let $\Phi$ be a monomial in $n=[r/4]$
variables. Let us denote by $W_\Phi$ its \emph{total weight},
defined by assigning weight $4i$ to its $i$-th variable.
We can use the monomial $\Phi$ to attach a real number to
the algebroid $A$.  Assume that $W_\Phi \le r$, and let $\Xi^{\g}$ 
be an equivariant twisted coycle of degree $r-W_\Phi$, i.e., $\Xi^{\g}\in\mathfrak Q^{r-W_\Phi}_G$.
Let $\Xi=\Xi^{\g}(0) \in \tilde C_A^{r-W_\Phi}$ be the zero degree term of $\Xi^{\g}$ as a function
of $\xi\in\g$;
it satisfies $\tilde\delta\Xi=0$.
We may   define
$$\Phi_\Xi(A) =(-2\pi)^{-m/2} \int_M \Phi(\lambda_2(A),\dots,\lambda_{2n}(A))\wedge  \Xi $$
where the $\lambda_i$ are the characteristic classes of the vector bundle $A$ as defined in
the previous section. This number only depends on the Lie algebroid $A$ and on the cohomology class $\Xi$. We also define
an element in $\sym^\bullet(\g^\ast)$
\begin{eqnarray}\label{eqpol}\Phi^{\g}_\Xi(A) &=&(-2\pi)^{-m/2} \int_M \Phi(\lambda_2^{\g}(A),\dots,\lambda^{\g}_{2n}(A))\wedge \Xi^{\g}  \nonumber
\\ &=& (-2\pi)^{-m/2} \int_M \Phi(\varsigma_2({\mathcal R}_\eta+\mu),\dots,\varsigma_{2n}({\mathcal R}_\eta+\mu))\wedge \Xi^{\g} \,. \end{eqnarray}
Of course, $\Phi_\Xi(A) = \Phi^{\g}_\Xi(A) (0)$.
 One has the following Bott-type result, which follows straightforwardly from Theorem
\ref{loc}.
\begin{thm} \label{bott1} Let $A$ be a rank $r$ Lie algebroid on an $m$-dimensional compact oriented manifold $M$, and  let $\alpha\in\Gamma(A)$ be any section
which  generates a circle action on $A$, and is such that $a(\alpha)$  has   isolated zeroes. Let $\Phi$ be a   
polynomial  in $n=[r/4]$ variables whose monomials have total weight $W_\Phi$.
Then if $r \ge m \ge W_\phi $ one has
\begin{equation}\label{polloc0}
\Phi_\Xi(A) = \sum_{x\in M_{a(\alpha)}}
\frac{\Phi(c_2(L_{a(\alpha),x}),\dots,c_{2n}(L_{a(\alpha),x}))\,p(\Xi^{\g})_{0}}{\operatorname{det}^{1/2} L_{a(\alpha),x}}
\end{equation}
where the classes $c_i(L_{a(\alpha),x})$ are the equivariant Chern classes of the endomorphism
$L_{a(\alpha),x}$ acting
on the tangent space $T_xM$ at a zero $x$ of $a(\alpha)$ (cf.~\cite{Bott67}).
\end{thm}
\begin{proof}  The r.h.s.~of Eq.~\eqref{polloc0} is computed from the r.h.s.~of the formula \eqref{locf}, taking into account two facts. First, we can evaluate the equivariant characteristic classes involved in $\Phi_{\Xi}^{\g}(A)$ by choosing in the principal bundle $GL(A)$ an equivariant $A$-connection induced by an ordinary equivariant connection $\zeta$ on the vector bundle $A$. In this way, we have
$$\Phi(\lambda^{\g}_2(A),\dots,\lambda^{\g}_{2n}(A))
= a^\ast(\Phi(\nu^{\g}_2(A),\dots,\nu^{\g}_{2n}(A))$$
where the classes $\nu_i^{\g}$ are equivariant Chern classes for the complexification of the
vector bundle $A$. Secondly, if $R_\zeta^{\g}$ is the equivariant curvature of $\zeta$, for every
symmetric elementary function $\varsigma_i$ and every zero $x$ of $a(\alpha)$, we have
$$\left(\varsigma_i(R^{\g}_\zeta) \right)_{0}(x)= c_i(L_{a(\alpha),x})\,.$$
\end{proof}

When $A=TM$, formula \eqref{polloc0}
 reduces to the ordinary Bott formula.

\section{Final remarks}\label{lastrem}
 As we shall show in a forthcoming paper \cite{BR},  Theorem \ref{bott1}   generalizes several localization formulas, for instance associated with the action of a holomorphic vector field on a complex manifold which is equivariantly lifted to an action
on a holomorphic vector bundle, see \cite{BB,Ch,CL77}, and reproduces in particular
Grothendieck's residue theorem.

On the other hand, our formula can be generalized in several directions. One of these
would be a localization formula for an equivariant cohomology associated with the action
of a Lie group on a Courant algebroid. This should encompass several formulas
recently appeared in the literature, mostly concerned with generalized Calabi-Yau
structures \cite{BCG07,Nitta06,Nitta07,Huribe06} and should reproduce our formula
when the Courant algebroid reduces to a Lie algebroid.

\end{document}